\documentclass[12pt]{article}
\usepackage{amsmath}
\usepackage{amssymb}
\usepackage{latexsym}
\usepackage{euscript}
\usepackage{amsfonts,amssymb}

\numberwithin{equation}{section}
\newtheorem{thm}[equation]{Theorem}

\newtheorem{lem}[equation]{Lemma}

\title{ Monogenic functions in the biharmonic boundary value problem}

\author{S.~V. Gryshchuk$^1$ \&\, S.~A. Plaksa$^2$\\
$^1$ $^,$ $^2$ Institute of Mathematics,\\ National Academy of
Sciences
of Ukraine,\\ Tereshchenkivska Str. 3, 01601, Kiev, Ukraine,\\
$^1$ serhii.gryshchuk@gmail.com, gryshchuk@imath.kiev.ua\\[2pt]
$^2$  plaksa@imath.kiev.ua}

\date{\ }

\begin{document}

\maketitle

\noindent {\bf Abstract:}  We consider a commutative algebra
$\mathbb{B}$ over the field of complex numbers with a basis
$\{e_1,e_2\}$ satisfying the conditions $(e_1^2+e_2^2)^2=0$,
$e_1^2+e_2^2\ne 0$. Let $D$ be a bounded domain in the Cartesian
plane $xOy$ and $D_{\zeta}=\{xe_1+ye_2 : (x,y)\in D\}$. Components
of every monogenic function
$\Phi(xe_1+ye_2)=U_{1}(x,y)\,e_1+U_{2}(x,y)\,ie_1+
U_{3}(x,y)\,e_2+U_{4}(x,y)\,ie_2$ having the classic derivative in
$D_{\zeta}$ are biharmonic functions in $D$, i.e.
$\Delta^{2}U_{j}(x,y)=0$ for $j=1,2,3,4$. We consider a
Schwarz-type boundary value problem for monogenic functions in a
simply connected domain $D_{\zeta}$. This problem is associated
with the following biharmonic problem: to find a biharmonic
function $V(x,y)$ in the domain $D$ when boundary values of its
partial derivatives $\partial V/\partial x$, $\partial V/\partial
y$ are given on the boundary $\partial D$. Using a hypercomplex
analog of the Cauchy type integral, we reduce the mentioned
Schwarz-type boundary value problem to a system of integral
equations on the real axes and establish sufficient conditions
under which this system has the Fredholm property.
\vspace{\baselineskip}


\section{Introduction}
\label{introd} Let $D$ be a bounded domain in the Cartesian plane
$xOy$, and let its boundary $\partial D$ be a closed smooth Jordan
curve. Let $\mathbb{R}$ be the set of real numbers.

Consider some boundary value problems for  {\em biharmonic
functions}\break $W \colon D\longrightarrow \mathbb{R}$ which have
continuous partial derivatives up to the fourth order inclusively
and satisfy the biharmonic equation in the domain $D$:
 \[  \Delta^{2}\,W(x,y)\equiv \frac{\partial^{4}W(x,y)}{\partial x^4}+
   2\frac{\partial^{4}W(x,y)}{\partial x^2\partial y^2}+
   \frac{\partial^4 W(x,y)}{\partial y^4} =0.\]

The {\em principal biharmonic problem} (cf., e.g., \cite[p.
194]{S3} and \cite[p. 13]{Mikhlin}) consists of finding a function
$W \colon \overline{D}\longrightarrow \mathbb{R}$ which is
continuous together with partial derivatives of the first order in
the closure $\overline{D}$ of the domain  $D$ and is biharmonic in
$D$, when its values and values of its outward normal derivative
are given on the boundary $\partial D$:
\begin{equation}\label{I_Big_pr}
W(x_0,y_0)=\omega_{1}(s), \quad \frac{\partial W}{\partial {\bf
n}}(x_0,y_0)= \omega_{2}(s)  \qquad \forall\, (x_0,y_0) \in
\partial D\,,
\end{equation}
 where $s$ is an arc coordinate of the point $(x_0,y_0)\in\partial D$.

In the case where $\omega_{1}$ is a continuously differentiable
function, the principal biharmonic problem is equivalent to the
following {\it biharmonic problem } (cf., e.g., \cite[p. 194]{S3}
and \cite[p. 13]{Mikhlin}) on finding a biharmonic function  $V :
D \longrightarrow\mathbb{R}$ with the following boundary
conditions:
\begin{equation}\label{osn_big_pr}
\begin{array}{l}
\displaystyle \lim\limits_{(x,y)\to(x_0,y_0),\, (x,y)\in D}
\frac{\partial V(x,y)}{\partial x}=
\omega_{3}(s)\,,\\[6mm]
\displaystyle \lim\limits_{(x,y)\to(x_0,y_0),\, (x,y)\in D}
\frac{\partial V(x,y)}{\partial y}= \omega_{4}(s)\qquad
\forall\, (x_0,y_0) \in \partial D\,,\\[6mm]
\displaystyle \int\limits_{\partial D}
\Bigl(\omega_{3}(s)\cos\angle({\bf
s},x)+\omega_{4}(s)\cos\angle({\bf s},y)\Bigr)\,ds=0\,.
\end{array}
\end{equation}
Here given boundary functions $\omega_{3}$, $\omega_{4}$ have
relations with given functions  $\omega_{1}$, $\omega_{2}$ of the
problem \eqref{I_Big_pr}, viz.,
\[
\omega_{3}(s)=\omega_{1}'(s)\cos\angle({\bf
s},x)+\omega_{2}(s)\cos\angle({\bf n},x)\,,  \]
\[
\omega_{4}(s)=\omega_{1}'(s)\cos\angle({\bf
s},y)+\omega_{2}(s)\cos\angle({\bf n},y)\,,  \] where ${\bf s}$
and ${\bf n}$ denote unit vectors of the tangent and the outward
normal to the boundary $\partial D$, respectively, and
$\angle(\cdot,\cdot)$ denotes an angle between an appropriate
vector ($\bf{s}$ or $\bf n$) and the positive direction of
coordinate axis ($x$ or $y$) indicated in the parenthesis.
Furthermore, solutions of the problems \eqref{I_Big_pr} and
\eqref{osn_big_pr} are related by the equality\,\,
$V(x,y)=W(x,y)+c$\,, where\, $c\in \mathbb{R}$\,.

A technique of using analytic functions of the complex variable
for solving the biharmonic problem is based on an expression of
biharmonic functions by the Goursat formula. This expression
allows to reduce the biharmonic problem to a certain boundary
value problem for a pair of analytic functions. Further,
expressing analytic functions via the Cauchy type integrals, one
can obtain a system of integro-differential equations in the
general case. In the case where the boundary $\partial D$ is a
Lyapunov curve, the mentioned system can be reduced to a system of
Fredholm equations. Such a scheme is developed (cf., e.g.,
\cite{Mush_upr,Lurie_engl,Mikh_int_eq_Th-El,Mikh_Dokl_43}) for
solving the main problems of the plane elasticity theory with
using a special biharmonic function which is called
the Airy stress function.

Another methods for reducing boundary value problems of the plane
elasticity theory to integral equations are developed in
\cite{Mikhlin,Mikh_kniga,Lu,Kakh,Kaland_79,Bogan}.

In this paper, for solving the biharmonic problem we develop a
method which is based on the relation between biharmonic functions
and monogenic functions taking values in a commutative algebra. We
use an expression of monogenic function by a hypercomplex analog
of the Cauchy type integral. Considering a Schwarz type boundary
value problem for monogenic functions that is associated with the
biharmonic problem, we develop a scheme of its reduction to a
system of Fredholm equations in the case where the boundary of
domain belongs to a class being wider than the class of Lyapunov
curves.
\section{Monogenic functions in a biharmonic algebra associated with the biharmonic equation}
\label{Mf_ass_bih_eq} V.~F. Kovalev and I.~P. Mel'nichenko
\cite{KM-BFf} considered an associative commutative
two-dimensional algebra $\mathbb B$ over the field of complex
numbers $\mathbb C$ with the following multiplication table for
basic elements $e_1, e_2$:
\begin{equation} \label{tab_umn_bb}
e_1^2=e_1,\quad e_2e_1=e_2,\quad e_2^2=e_1+2ie_2,
\end{equation}
where\, $i$\, is the imaginary complex unit. Elements $e_1, e_2$
satisfy the relations
\begin{equation} \label{tab_umn_ba}
(e_1^2+e_2^2)^2=0,\qquad e_1^2+e_2^2\ne 0,
\end{equation}
thereby, all functions $\Phi(\zeta)$ of the variable
$\zeta=x\,e_1+y\,e_2$, which have continuous derivatives up to the
fourth order inclusively, satisfy the equalities
\begin{equation} \label{mon_ud_bu}
\Delta^{2}\,\Phi(\zeta)=\Phi^{(4)}(\zeta)\,(e_1^2+e_2^2)^2=0\,.
\end{equation}
Therefore, components  $U_{j}\colon D\longrightarrow \mathbb{R}$,
$j=\overline{1,4}$, of the expression
\begin{equation}\label{mon-funk}
\Phi(\zeta)=U_{1}(x,y)\,e_1+U_{2}(x,y)\,ie_1+
U_{3}(x,y)\,e_2+U_{4}(x,y)\,ie_2
\end{equation}
are biharmonic functions.

Similarly to \cite{KM-BFf}, by the {\em biharmonic plane} we call
a linear span $\mu:=\{\zeta=x\,e_1+y\,e_2: x,y \in \mathbb{R}\}$
of the elements $e_1$, $e_2$ satisfying (\ref{tab_umn_bb}). With a
domain $D$ of the Cartesian plane $xOy$ we associate the congruent
domain $D_{\zeta}:= \{\zeta=xe_1+ye_2 : (x,y)\in D\}$ in the
biharmonic plane $\mu$.


We say that a function $\Phi\colon
D_{\zeta}\longrightarrow\mathbb{B}$ is {\em monogenic}\/ in a
domain $D_{\zeta}$ if at every point $\zeta\in D_{\zeta}$
 there exists the derivative of the function $\Phi$:
\[\Phi'(\zeta):=\lim\limits_{h\to 0,\, h\in \mu}
\bigl(\Phi(\zeta+h)-\Phi(\zeta)\bigr)\,h^{-1}\,.\]

It is proved in \cite{KM-BFf} that a function  $\Phi\colon
D_{\zeta}\longrightarrow \mathbb{B}$ is monogenic in a domain
$D_{\zeta}$ if and only if components  $U_{j}\colon
D\longrightarrow \mathbb{R}$, $j=\overline{1,4}$, of the
expression \eqref{mon-funk} are differentiable in the domain $D$
and the following analog of the Cauchy -- Riemann conditions  is
satisfied:
\begin{equation}\label{usl_K_R}
\frac{\partial \Phi(\zeta)}{\partial y}=\frac{\partial
\Phi(\zeta)}{\partial x}\,e_2\quad \forall\, \zeta=xe_1+e_2y \in
D_{\zeta}\,.
\end{equation}

It is established in \cite{GrPl_umz-09} that every monogenic
function $\Phi \colon D_{\zeta}\longrightarrow \mathbb{B}$ has
derivatives $\Phi^{(n)}(\zeta)$ of all orders $n$ in the domain
$D_{\zeta}$ and, therefore, it satisfies the equalities
\eqref{mon_ud_bu}. At the same time, every biharmonic in  $D$
function $U(x,y)$ is the first component $U_{1}\equiv U$ in the
expression (\ref{mon-funk}) of a certain monogenic function  $\Phi
: D_{\zeta}\longrightarrow \mathbb{B}$ and, moreover, all such
functions $\Phi$ are found in \cite{GrPl_umz-09} in an explicit
form.
\section{Boundary value problem for monogenic functions that is associated with the biharmonic problem}
\label{BVP_ass_Bih-Pr} Let $\Phi_{1}$ be monogenic in $D_{\zeta}$
function having the sought-for function
 $V(x,y)$ of the problem \eqref{osn_big_pr} as the first component:
 \[\Phi_{1}(\zeta)=V(x,y)\,e_1+V_{2}(x,y)\,ie_1+V_{3}(x,y)\,e_2+V_{4}(x,y)\,ie_2\,.\]
 It follows from the condition  (\ref{usl_K_R}) for $\Phi=\Phi_1$
that $\partial V_{3}(x,y)/\partial x=\partial V(x,y)/\partial y$.
Therefore,
\begin{equation}\label{osn-(1-3)}
\Phi_{1}'(\zeta)=\frac{\partial V (x,y)}{\partial x}\,e_1+
\frac{\partial V_{2}(x,y)}{\partial x}\,ie_1+ \frac{\partial
V(x,y)}{\partial y}\,e_2+ \frac{\partial V_{4}(x,y)}{\partial
x}\,ie_2\,
\end{equation}
 and, as consequence, we conclude that the biharmonic problem with boundary conditions (\ref{osn_big_pr})
is reduced to the boundary value problem on finding a monogenic in
$D_{\zeta}$ function $\Phi\equiv\Phi_{1}'$ when values of two
components $U_1=\partial V(x,y)/\partial x$ and $U_3=\partial
V(x,y)/\partial y$ of the expression (\ref{mon-funk}) are given on
the boundary  $\partial D_{\zeta}$ of the domain  $D_{\zeta}$.

As in \cite{IJPAM_13}, by the (1-3)-{\it problem} we shall call
the problem on finding a monogenic function  $\Phi \colon
D_{\zeta}\longrightarrow \mathbb{B}$ when values of components
$U_1$, $U_3$ of the expression (\ref{mon-funk}) are given on the
boundary $\partial D_{\zeta}$:
\[U_{1}(x_0,y_0)=u_{1}(\zeta_0)\,,\quad
U_{3}(x_0,y_0)=u_{3}(\zeta_0) \qquad
\forall\,\zeta_0:=x_0e_1+y_0e_2\in\partial D_{\zeta},\]
 where
$u_{1}(\zeta_0)\equiv\omega_3(s)$\,  and \,
$u_{3}(\zeta_0)\equiv\omega_4(s)$.

Problems of such a type on finding a monogenic function with given
boundary values of two its components were posed by V.~F. Kovalev
\cite{Kovalov} who called them by {\it biharmonic Schwarz
problems}, because their formulations are analogous in a certain
sense to the classic Schwarz problem on finding an analytic
function of the complex variable when values of its real part are
given on the boundary of domain. V.~F. Kovalev \cite{Kovalov}
stated a sketch of reduction of biharmonic Schwarz problems to
integro-differential equations with using conformal mappings and
expressions of monogenic functions via analytic functions of the
complex variable.

In  \cite{IJPAM_13}, we investigated the (1-3)-problem for cases
where $D_{\zeta}$ is either an upper half-plane or a unit disk in
the biharmonic plane. Its solutions were found in explicit forms
with using of some integrals analogous to the classic Schwarz
integral.

In \cite{G_conf_14}, a certain scheme was proposed for reducing
the (1-3)-problem in a simply connected domain with sufficiently
smooth boundary to a suitable boundary value problem in a disk
with using power series and conformal mappings in the complex
plane.

Hypercomplex methods for investigating the biharmonic equation
were developed in the papers
\cite{Sodbero,Vignaux_uno,Edenhofer,Snyder_Cont_82,Snyder_comm_alg_82,Bock_Gurl}
also.
\section{Biharmonic Cauchy type integral}
\label{bihCint}
 Consider the biharmonic Cauchy type integral
\begin{equation}\label{int_t_Cauchy}
\Phi(\zeta)=\frac{1}{2\pi i}\int\limits_{\partial
D_\zeta}\varphi(\tau)(\tau-\zeta)^{-1}d\tau
\end{equation}
with a continuous density $\varphi\colon\partial
D_\zeta\longrightarrow\mathbb{B}$. The integral
(\ref{int_t_Cauchy}) is a monogenic function in both domains
$D_{\zeta}$ and $\mu\setminus\overline{D_{\zeta}}$.

We use the euclidian norm $\|a\|:=\sqrt{|z_1|^2+|z_2|^2}$ in the
algebra $\mathbb{B}$, where\, $a=z_1e_1+z_2e_2$ and
 $z_1, z_2\in \mathbb{C}$.
We use also the modulus of continuity of a function $\varphi
\colon \partial D_\zeta\longrightarrow~\mathbb{B}$:
\[\omega(\varphi,\varepsilon):=\sup\limits_{\tau_1,\tau_2\in\partial
D_\zeta,\,
\|\tau_1-\tau_2\|\le\varepsilon}\|\varphi(\tau_1)-\varphi(\tau_2)\|\,.\]
Consider a singular integral which is understood in the sense of
its Cauchy principal value:
\[\int\limits_{\partial D_\zeta}\varphi(\tau)(\tau-\zeta_0)^{-1}d\tau:=\lim\limits_{\varepsilon\rightarrow0}\int\limits_
{\partial D_\zeta\setminus\partial D_\zeta^{\varepsilon}
(\zeta_0)}\varphi(\tau)(\tau-\zeta_0)^{-1}d\tau,\]
where $\zeta_0\in\partial D_\zeta$, 
$\partial D_\zeta^{\varepsilon}(\zeta_0):= \{\tau\in\partial
D_\zeta:\|\tau-\zeta_0\|\leq\varepsilon\}$.

The following theorem can be proved in a similar way as an
appropriate theorem in the complex plane (cf., e.g.,
\cite{Gah,Magnar_For_Sokh}). It presents sufficient conditions for
the existence of limiting values
\[\Phi^{+}(\zeta_0):=\lim\limits_{\zeta\to\zeta_0,\, \zeta \in
D_{\zeta}}\Phi(\zeta), \qquad
\Phi^{-}(\zeta_0):=\lim\limits_{\zeta\to\zeta_0,\, \zeta\in
\mu\setminus\overline{D_{\zeta}}}\Phi(\zeta) \]
 of the biharmonic Cauchy type integral in any point $\zeta_0\in\partial
 D_{\zeta}$.

\begin{thm}
\label{form_Sokh_D} Let the modulus of continuity of a function
$\varphi \colon \partial D_{\zeta} \longrightarrow \mathbb{B}$
satisfy the Dini condition
\begin{equation} \label{usl_Dini}
\int\limits_{0}^{1}\frac{\omega(\varphi,\eta)}{\eta}\,d\,\eta<\infty\,.
\end{equation}
Then the integral  \eqref{int_t_Cauchy} has limiting values
$\Phi^{\pm}(\zeta_0)$ in any point $\zeta_0 \in \partial
D_{\zeta}$ that are represented by the Sokhotski--Plemelj
formulas: \begin{equation} \label{int_t_Cauchy+}
\begin{array}{l}
\displaystyle \Phi^{+}(\zeta_0)=\frac{1}{2}\,\varphi(\zeta_0)+
 \frac{1}{2\pi i}\int\limits_{\partial D_{\zeta}}\varphi(\tau)(\tau-\zeta_0)^{-1}\,d\tau\,,\\[10mm]
\displaystyle \Phi^{-}(\zeta_0)=-\frac{1}{2}\,\varphi(\zeta_0)+
 \frac{1}{2\pi i}\int\limits_{\partial D_{\zeta}}\varphi(\tau)(\tau-\zeta_0)^{-1}\,d\tau\,.
  \end{array}
\end{equation}
\end{thm}

\section{Scheme for reducing the (1-3)-problem to a system of integral equations}
\label{Reduct13} Let the functions  $u_1 \colon \partial
D_{\zeta}\longrightarrow \mathbb{R}$, $u_3 \colon  \partial
D_{\zeta}\longrightarrow \mathbb{R}$ satisfy conditions of the
type \eqref{usl_Dini}.

We shall find solutions of the (1-3)-problem in the class of
functions represented in the form
\begin{equation} \label{int_t_Cauchy-1-3}
 \Phi(\zeta)= \frac{1}{2\pi
i}\int\limits_{\partial D_{\zeta}}\Bigl(\varphi_{1}(\tau)e_1+
\varphi_{3}(\tau)e_2\Bigr)(\tau-\zeta)^{-1}\,d\tau \qquad
\forall\, \zeta \in D_{\zeta},
\end{equation}
where the functions $\varphi_1 : \partial D_{\zeta}\rightarrow
\mathbb{R}$ and $\varphi_3 : \partial D_{\zeta}\rightarrow
\mathbb{R}$ satisfy conditions of the type \eqref{usl_Dini}.

Then, by Theorem~\ref{form_Sokh_D},  the following equality is
valid for any  $\zeta_0 \in\partial D_{\zeta}$:
 \begin{equation} \label{int_t_Cauchy-13}
 \Phi^{+}(\zeta_0):= \frac{1}{2}\Bigl(\varphi_{1}(\zeta_0)e_1+\varphi_{3}(\zeta_0)e_2\Bigr)+ \frac{1}{2\pi
i}\int\limits_{\partial D_{\zeta}}\Bigl(\varphi_{1}(\tau)e_1+
\varphi_{3}(\tau)e_2\Bigr)(\tau-\zeta_0)^{-1}\,d\tau\,. 
\end{equation}

By $D_{z}$ we denote the domain in $\mathbb{C}$ which is congruent
to the domain $D$, i.e. $D_{z}:=\{z=x+iy\in\mathbb{C}:(x,y)\in D
\}$. We shall use a conformal mapping  $z=\tau(t)$ of the upper
half-plane  $\{t\in \mathbb{C}: \mathrm{Im}\,t>0\}$ onto the
domain $D_{z}$. Denote $\tau_1(t):=\mathrm{Re}\,\tau(t)$,
$\tau_2(t):=\mathrm{Im}\,\tau(t)$.

Inasmuch as the mentioned conformal mapping is continued to a
homeomorphism between the closures of corresponding domains, the
function
\begin{equation}\label{tilde_tay}
\widetilde{\tau}(s):=\tau_{1}(s)e_1+\tau_{2}(s)e_2 \qquad
\forall\, s \in \overline{\mathbb{R}}
\end{equation}
 generates a homeomorphic mapping of the extended real axis  $\overline{\mathbb{R}}:=\mathbb{R}\cup\{\infty\}$ onto the curve
$\partial D_{\zeta}$.

Introducing the function
\begin{equation} \label{fun_g}
 g(s) :=g_1(s)e_1+g_{3}(s)e_2 \qquad  \forall\, s\in \overline{\mathbb{R}}\,,
\end{equation}
where \, $g_{j}(s):=\varphi_{j}\left( \widetilde{\tau}(s)\right)$
for $j\in\{1,3\}$, we rewrite the equality (\ref{int_t_Cauchy-13})
in the form
 \begin{equation} \label{C1-3}
\Phi^+(\zeta_0)=\frac{1}{2}\,g(t)+\frac{1}{2\pi
i}\int\limits_{-\infty}^{\infty} g(s)
   \bigl(\widetilde{\tau}(s)-\widetilde{\tau}(t)\bigr)^{-1}\,\,{\widetilde{\tau}}\,'(s)\,ds
   \qquad \forall\,t\in\mathbb{R}\,,
 \end{equation}
where the integral is understood in the sense of its Cauchy
principal value (cf., e.g., \cite{Gah}) and a correspondence
between the points $\zeta_0\in\partial
D_{\zeta}\setminus\{\widetilde{\tau}(\infty)\}$ and
$t\in\mathbb{R}$ is given by the equality
 $\zeta_0=\widetilde{\tau}(t)$.

To transform the expression under an integral sign in the equality
\eqref{C1-3} we use the equalities
\[\bigl(\widetilde{\tau}(s)-\widetilde{\tau}(t)\bigr)^{-1}= \frac{1}{\tau(s)-\tau(t)}+
\frac{i\bigl(\tau_{2}(s)-\tau_{2}(t)\bigr)}{2\bigl(\tau(s)-\tau(t)\bigr)^{2}}\,\rho\,,\]
\[{\widetilde{\tau}}\,'(s)={\tau}'(s)-\frac{i\tau_{2}'(s)}{2}\,\rho\,,\]
where
\begin{equation} \label{rho}
\rho:=2e_1+2ie_2
\end{equation}
is a nilpotent element of the algebra $\mathbb{B}$ because
$\rho^2=0$. Thus, we transform the equality (\ref{C1-3}) into the
form
\[\Phi^+(\zeta_0)=\frac{1}{2}\,g(t)+ \frac{1}{2\pi
i}\int\limits_{-\infty}^{\infty}g(s)k(t,s)\,ds+ \frac{1}{2\pi
i}\int\limits_{-\infty}^{\infty}
g(s)\frac{1+st}{(s-t)(s^2+1)}\,ds\,,\]
 where
\begin{equation} \label{k}
k(t, s) = k_{1}(t, s)e_1+i\rho \, k_{2}(t, s)\,,
\end{equation}
\begin{equation} \label{k1}
k_{1}(t,s):=\frac{\tau'(s)}{\tau(s)-\tau(t)}-\frac{1+st}{(s-t)(s^2+1)}\,,
\end{equation}
\begin{equation} \label{k2}
k_{2}(t,
s):=\frac{\tau'(s)\bigl(\tau_{2}(s)-\tau_{2}(t)\bigr)}{2\bigl(\tau(s)-\tau(t)\bigr)^{2}}-
\frac{\tau_{2}'(s)}{2\bigl(\tau(s)-\tau(t)\bigr)}\,.
\end{equation}

We use the notations $U_{j}\left[a\right]:=a_{j}$,
$j=\overline{1,4}$, where $a_{j}\in\mathbb{R}$ is the coefficient
in the decomposition of element $a=a_1e_1+a_2 ie_1+a_3e_2+a_4
ie_2\in\mathbb{B}$ with respect to the basis $\{e_1,e_2\}$.

In order to single out $U_{1}\left[\Phi^+(\zeta_0)\right]$,
$U_{3}\left[\Phi^+(\zeta_0)\right]$ we use the equalities
\eqref{fun_g}, \eqref{rho}, \eqref{k} and get the decomposition of
the following expression with respect to the basis $\{e_1,e_2\}$:
\[g(s)k(t,s)=\bigl(g_1(s)e_1+g_{3}(s)e_2\bigr)\bigl(k_{1}(t,s)e_1+i(2e_1+2ie_2)\, k_{2}(t,s)\bigr)=\]
\[=\Bigl(g_1(s)\bigl(k_{1}(t,s)+2ik_{2}(t,s)\bigr)-2g_{3}(s)k_{2}(t,s)\Bigr)e_1+\]
 \[+\Bigl(g_3(s)\bigl(k_{1}(t,s)-2ik_{2}(t,s)\bigr)-2g_{1}(s)k_{2}(t,s)\Bigr)e_2\,.\]

Now, we single out $U_{1}\left[\Phi^+(\zeta_0)\right]$,
$U_{3}\left[\Phi^+(\zeta_0)\right]$ and obtain the following
system of integral equations for finding the functions $g_1$ and
$g_3$:
\begin{equation} \label{syst-int-ur}
\begin{array}{l}
\displaystyle
 U_{1}\left[\Phi^+(\zeta_0)\right]\equiv\frac{1}{2}\,g_{1}(t)+\frac{1}{2\pi}\int\limits_{-\infty}^{\infty}
g_{1}(s)\Bigl(\mathrm{Im}\,k_1(t,s)+2\mathrm{Re}\,k_2(t,s)\Bigr)\,ds-\\[4mm]
 \hspace*{27mm} \displaystyle -\frac{1}{\pi}\int\limits_{-\infty}^{\infty}g_{3}(s)\mathrm{Im}\,k_2(t,s)\,ds=\widetilde{u}_{1}(t),\\[4mm]
 \displaystyle U_{3}\left[\Phi^+(\zeta_0)\right]\equiv\frac{1}{2}\,g_{3}(t)-\frac{1}{\pi}\int\limits_{-\infty}^{\infty}g_{1}(s)\mathrm{Im}\,k_2(t,s)\,ds+\\[4mm]
 \hspace*{17mm}
 \displaystyle +\frac{1}{2\pi}\int\limits_{-\infty}^{\infty}
g_{3}(s)\Bigl(\mathrm{Im}\,k_1(t,s)-2\mathrm{Re}\,k_2(t,s)\Bigr)\,ds=\widetilde{u}_{3}(t)
\qquad \forall\,t\in\mathbb{R}\,,  \end{array}
\end{equation}
where
$\widetilde{u}_{j}(t):=u_{j}\bigl(\widetilde{\tau}(t)\bigr)$,\,\,
$j\in\{1,3\}$.

Below, we shall state conditions which are sufficient for
compactness of integral operators on the left-hand sides of
equations of the system \eqref{syst-int-ur}.
\section{Auxiliary Statements}\label{aux_stat}
For a function $\varphi : \gamma \longrightarrow\mathbb{C}$ which
is continuous on the curve $\gamma\subset\mathbb{C}$,  a modulus
of continuity is defined by the equality
\[
\omega_{\gamma}(\varphi,\varepsilon):=\sup_{t_1, t_2 \in \gamma,\,
|t_1-t_2| \le \varepsilon}|\varphi(t_1)-\varphi(t_2)|\,.
\]

Consider the conformal mapping $\sigma(T)$ of the unit disk
$\{T\in\mathbb{C}: |T|<1\}$ onto the domain $D_{z}$ such that
$\tau(t)=\sigma\left(\frac{t-i}{t+i}\right)$ for all  $t\in\{t\in
\mathbb{C}: \mathrm{Im}\,t>0\}$. Denote
$\sigma_1(T):=\mathrm{Re}\,\sigma(T)$,
$\sigma_2(T):=\mathrm{Im}\,\sigma(T)$.

Assume that the conformal mapping $\sigma(T)$ has the continuous
contour derivative $\sigma'(T)$ on the unit circle
$\Gamma:=\{T\in\mathbb{C}: |T|=1\}$ and $\sigma'(T)\ne 0$ for all
$T\in\Gamma$. Then there exist constants $c_1$ and $c_2$ such that
the following inequalities are valid:
\begin{equation} \label{otdel_0}
0<c_{1}\le \left|\frac{\sigma(S)-\sigma(T)}{S-T}\right| \le
c_{2}\,.
\end{equation}
In this case, for all  $S, T_1, T_2\in\Gamma$ such that
$0<|S-T_1|<|S-T_2|$ the following estimates are also valid:
\begin{equation} \label{prir_d_1_k}
  \left| \frac{\sigma_{j}(S)-\sigma_{j}(T_1)}{S-T_1} - \frac{\sigma_{j}(S)-\sigma_{j}(T_2)}{S-T_2}
  \right| \le c\,\frac{\omega_{\Gamma}(\sigma', |S-T_2|)}{|S-T_2|} \, |T_1-T_2|\,, \qquad
  j=1,2,
\end{equation}
where the constant $c$ does not depend on  $S$, $T_1$, $T_2$.
Estimates of the same type for the function $\sigma(T)$ are
corollaries of the estimates \eqref{prir_d_1_k}. They are adduced
in \cite{P5}.

It follows from  \eqref{prir_d_1_k} that the inequalities
  \begin{equation} \label{sigma'_1}
 \left|\sigma'_j(S)- \frac{\sigma_{j}(S)-\sigma_{j}(T)}{S-T}  \right| \le c \, \omega_{\Gamma}(\sigma', |S-T|)\,,\qquad
 j=1,2,
  \end{equation}
are fulfilled for all  $S, T \in \Gamma$, $S \ne T$, where the
constant \, $c$\, does not depend on $S$ and $T$. Then an
inequality of the same type for the function $\sigma(T)$ is
certainly fulfilled.

Let $C(\overline{\mathbb{R}}\,)$ denote the Banach space of
functions $g \colon
\overline{\mathbb{R}}\longrightarrow\mathbb{C}$ that are
continuous on the extended real axis $\overline{\mathbb{R}}$ with
the norm
$\|g\|_{C(\overline{\mathbb{R}}\,)}:=\sup\limits_{t\in\mathbb{R}}|g(t)|$.

\begin{lem} \label{dini_K1}  Let $g\in C(\,\overline{\mathbb{R}}\,)$ and the conformal mapping
$\sigma(T)$ have the nonvanishing continuous contour derivative
$\sigma'(T)$ on the circle $\Gamma$, and its modulus of continuity
satisfy the Dini condition
 \begin{equation} \label{Dini_circle}
\int\limits_{0}^{1}\frac{\omega_{\Gamma}(\sigma',
\eta)}{\eta}\,d\eta<\infty.
\end{equation}
Then for $0<\varepsilon<1/4$ and  $t\in\mathbb{R}$ the following
estimates are true:
 \[\left|\,\int\limits_{-\infty}^{\infty}g(s)\,k_{j}(t+\varepsilon,s)\,ds-\int\limits_{-\infty}^{\infty}g(s)\,k_{j}(t,s)
 \,ds\,\right| \le\]
\begin{equation} \label{ner_K1'}
 c\,\|g\|_{C(\overline{\mathbb{R}}\,)}
 \,\epsilon\int\limits_{0}^{2}\frac{\omega_{\Gamma}(\sigma', \eta)}{\eta(\eta+\epsilon)}\,
 d\eta\,,\qquad j=1,2\,,
      \end{equation}
where \,\, $\epsilon:=\varepsilon / (t^2+1)$ and the constant\,
$c$ does not depend on $t$ and $\varepsilon$.
\end{lem}
{\bf Proof.} Let us prove the inequality (\ref{ner_K1'}) for
$j=1$. Denote $S:=\frac{s-i}{s+i}$\,, $T:=\frac{t-i}{t+i}$\,,
$d(S,T):=\frac{\sigma(S)-\sigma(T)}{S-T}$\,. Taking into account
the equalities
\begin{equation} \label{S-T}
S-T =\frac{2i(s-t)}{(s+i)(t+i)}\,,
\end{equation}
\[\tau'(s)=\frac{2i}{(s+i)^{2}}\,\sigma'(S)\,,\] we transform the expression
\[\frac{\tau'(s)}{\tau(s)-\tau(t)}=\frac{\frac{2i}{(s+i)^{2}}\,\sigma'(S)}{(S-T)d(S,T)}=
\left(\frac{1+st}{(s-t)(s^2+1)}+i\,\frac{1}{s^2+1}\right)\frac{\sigma'(S)}{d(S,T)}\]
and represent the function  $k_{1}(t,s)$ in the form
\begin{equation} \label{k1-s-m}
 k_{1}(t,s)= m_{1}(t,s)+i\,m_{2}(t,s),
\end{equation}
where
\[ m_{1}(t,s):=\frac{1+st}{(s-t)(s^2+1)}\frac{\sigma'(S)-d(S,T)}{d(S,T)}\,, \quad
m_{2}(t,s):=\frac{\sigma'(S)}{d(S,T)(s^2+1)}\,.
\]

Further, representing the integral
\[I[g,k_1](t):=
\int\limits_{-\infty}^{\infty}g(s)k_1(t,s)\,ds\qquad \forall\, t
\in \mathbb{R}\]
by the sum of two integrals
\begin{equation} \label{K1-int}
I[g,k_1](t)=I[g,m_1](t)+iI[g,m_2](t)\qquad  \forall\, t \in
\mathbb{R},
\end{equation}
we shall obtain estimates of the type (\ref{ner_K1'}) for each of
integrals $I[g,m_1](t)$, $I[g,m_2](t)$.

For the integral  $I[g,m_1](t)$, we have
\[
\Bigl|I[g,m_1](t+\varepsilon)- I[g,m_1](t)\Bigr| \le
\int\limits_{t-2\varepsilon}^{t+2\varepsilon}\left|g(s)\right|\left|m_{1}(t,s)\right|\,
ds+\]
\[+\int\limits_{t-2\varepsilon}^{t+2\varepsilon}\left|g(s)\right||m_{1}(t+\varepsilon,s)|\, ds+
\]
\[+\left(\int\limits_{-\infty}^{t-2\varepsilon}+\int\limits_{t+2\varepsilon}^{\infty}\right)
\left|g(s)\right|\bigl|m_{1}(t,s)-m_{1}(t+\varepsilon,s)\bigr|\,
ds=:J_1+J_2+J_3.\]

Taking into account the relations \eqref{otdel_0},
\eqref{sigma'_1}, \eqref{S-T},  we obtain
\[J_{1}\le c\,\|g\|_{C(\overline{\mathbb{R}}\,)}  \int\limits_{t-2\varepsilon}^{t+2\varepsilon} \frac{\omega_{\Gamma}(\sigma', |S-T|)}{|S-T|} 
\frac{\left(1+|s||t|\right)}{\sqrt{s^2+1}\sqrt{t^2+1}} \,
\frac{ds}{s^2+1} \le \]
\[ \le c \,\|g\|_{C(\overline{\mathbb{R}}\,)} \int\limits_{t-2\varepsilon}^{t+2\varepsilon}
\frac{\omega_{\Gamma}(\sigma', |S-T|)}{|S-T|} \,
\frac{ds}{s^2+1}\le\]
\[ \le c\,\|g\|_{C(\overline{\mathbb{R}}\,)} \int\limits_{0}^{8\epsilon} \frac{\omega_{\Gamma}(\sigma',
\eta)}{ \eta}\, d\eta \le c\,\|g\|_{C(\overline{\mathbb{R}}\,)}
 \,\epsilon\int\limits_{0}^{2}\frac{\omega_{\Gamma}(\sigma', \eta)}{\eta(\eta+\epsilon)}\,
 d\eta. \]
Here and below in the proof, by\, $c$\, we denote constants whose
values are independent of  $t$ and $\varepsilon$, but, generally
speaking, may be different even within a single chain of
inequalities.

The integral $J_{2}$ is similarly estimated.

To estimate the integral  $J_{3}$, take into consideration the
point $T_1:=\frac{t+\varepsilon-i}{t+\varepsilon+i}$. Using the
equality
\[m_{1}(t,s)-m_{1}(t+\varepsilon,s)=\frac{-\varepsilon}{(s-t)(s-t-\varepsilon)}\,\frac{\sigma'(S)-d(S,T)}{d(S,T)}\,+\]
\[+\frac{1+s(t+\varepsilon)}{(s-t-\varepsilon)(s^2+1)}\,\frac{\bigl(\sigma'(S)-d(S,T)\bigr)\bigl(d(S,T_1)-d(S,T)\bigr)}{d(S,T)d(S,T_1)}\,+\]
\[+\frac{1+s(t+\varepsilon)}{(s-t-\varepsilon)(s^2+1)}\,\frac{d(S,T_1)-d(S,T)}{d(S,T_1)}\,,\]
we estimate  $J_3$ by the sum of three integrals:
\[J_3\le \varepsilon\left(\int\limits_{-\infty}^{t-2\varepsilon}+\int\limits_{t+2\varepsilon}^{\infty}\right)
\left|g(s)\right|\frac{\bigl|\sigma'(S)-d(S,T)\bigr|}{|\,d(S,T)||s-t||s-t-\varepsilon|}\,ds\,+
\]
\[
+\left(\int\limits_{-\infty}^{t-2\varepsilon}+\int\limits_{t+2\varepsilon}^{\infty}\right)
\left|g(s)\right|\left|\,d(S,T_1)-d(S,T)\right|\times\]
\[\times\frac{\bigl|\sigma'(S)-d(S,T)\bigr|\bigl(1+|\,s||\,t+\varepsilon|\bigr)}{|\,d(S,T)|
|\,d(S,T_1)||s-t-\varepsilon|}\, \frac{ds}{s^2+1}\,+
\]
\[+\left(\int\limits_{-\infty}^{t-2\varepsilon}+\int\limits_{t+2\varepsilon}^{\infty}\right)
\left|g(s)\right|\frac{\left|\,d(S,T_1)-d(S,T)\right|
\bigl(1+|s||\,t+\varepsilon|\bigr)}{|\,d(S,T_1)||s-t-\varepsilon|}\,
\frac{ds}{s^2+1}=:\]
\[=:\sum\limits_{j=1}^{3}J_{3,\,j}\,.\]

Taking into account the inequalities  $|s-t|\le 2|s-t-\varepsilon|
\le 3|s-t|$ for all  $s\in(-\infty, t-2\varepsilon] \cup
[t+2\varepsilon, +\infty)$ and the relations (\ref{otdel_0}),
\eqref{sigma'_1}, (\ref{S-T}), we obtain
\[J_{3,1}\le c \,\|g\|_{C(\overline{\mathbb{R}}\,)}\, \varepsilon\left(\int\limits_{-\infty}^{t-2\varepsilon}+\int\limits_{t+2\varepsilon}^{\infty}\right)
\frac{\omega_{\Gamma}(\sigma', |S-T|)}{|s-t|^{2}}\,ds =
\]
\[
=c\,\|g\|_{C(\overline{\mathbb{R}}\,)}\,\,\frac{\varepsilon}{t^2+1}
\left(\int\limits_{-\infty}^{t-2\varepsilon}+\int\limits_{t+2\varepsilon}^{\infty}\right)
\frac{\omega_{\Gamma}(\sigma', |S-T|)}{|S-T|^{2}} \,
\frac{ds}{s^2+1} \le
\]
\[
\le c\,\|g\|_{C(\overline{\mathbb{R}}\,)}\, \epsilon
\int\limits_{\epsilon}^{2}\frac{\omega_{\Gamma}(\sigma',\eta)}{\eta^2}\,
d\eta \le c\,\|g\|_{C(\overline{\mathbb{R}}\,)}
\,\,\epsilon\int\limits_{0}^{2}\frac{\omega_{\Gamma}(\sigma',
\eta)}{\eta(\eta+\epsilon)}\,
 d\eta\,.\]

Using the inequalities  (\ref{prir_d_1_k}) and properties of a
modulus of continuity (cf., e.g., \cite[p.~176]{Dzjadyk}), we
obtain the following inequalities similarly to the estimation of
$J_{3,1}$:
\[J_{3,2} \le c\,\|g\|_{C(\overline{\mathbb{R}}\,)}\left(\int\limits_{-\infty}^{t-2\varepsilon}+\int\limits_{t+2\varepsilon}^{\infty}\right)\,
\frac{\left|d(S,T_1)-d(S,T)\right|\bigl(1+|\,s||\,t+\varepsilon|\bigr)}{|s-t-\varepsilon|}
\, \frac{ds}{s^2+1}\le \]
\[\le c\,\|g\|_{C(\overline{\mathbb{R}}\,)}  |T_1-T|\left(\int\limits_{-\infty}^{t-2\varepsilon}+\int\limits_{t+2\varepsilon}^{\infty}\right)\,
\frac{\omega(\sigma', |S-T|)}{ |S-T|}\,
\frac{\bigl(1+|\,s||\,t+\varepsilon|\bigr)}{|s-t|} \,
\frac{ds}{s^2+1} \le \]
\[\le c\,\|g\|_{C(\overline{\mathbb{R}}\,)}  |T_1-T|\left(\int\limits_{-\infty}^{t-2\varepsilon}+\int\limits_{t+2\varepsilon}^{\infty}\right)\,
\frac{\omega(\sigma', |S-T|)}{|S-T|^2}
\frac{\bigl(1+|\,s||\,t+\varepsilon|\bigr)}{\sqrt{s^2+1}\sqrt{t_1^2+1}}
\, \frac{ds}{s^2+1} \le \]
\[ \le c\,\|g\|_{C(\overline{\mathbb{R}}\,)}\, \epsilon \left(\int\limits_{-\infty}^{t-2\varepsilon}+\int\limits_{t+2\varepsilon}^{\infty}\right)\,
 \frac{\omega(\sigma', |S-T|)}{ |S-T|^{2}} \, \frac{ds}{s^2+1} \le \]
 \[ \le  c\,\|g\|_{C(\overline{\mathbb{R}}\,)}
\,\,\epsilon\int\limits_{0}^{2}\frac{\omega_{\Gamma}(\sigma',
\eta)}{\eta(\eta+\epsilon)}\,
 d\eta\,. \]

The integral $J_{3,3}$ is similarly estimated. An estimate of the
type (\ref{ner_K1'}) for $I[g,m_1](t)$ follows from the obtained
estimates.

By the scheme used above for estimating the integral
$I[g,m_1](t)$, we get the estimate
\[\Bigl|I[g,m_2](t+\varepsilon)- I[g,m_2](t)\Bigr| \le
c\,\|g\|_{C(\overline{\mathbb{R}}\,)}\,\, \epsilon\,, \]
 whence an
estimate of the type (\ref{ner_K1'}) follows for $I[g,m_2](t)$.
Thus, the inequality (\ref{ner_K1'}) is proved for $j=1$.

Let us prove the inequality (\ref{ner_K1'}) for $j=2$. Denote
$d_{j}(S,T):=\frac{\sigma_{j}(S)-\sigma_{j}(T)}{S-T}$ for $j=1,2$.

Similarly to the expression \eqref{k1-s-m}, we represent the
function  $k_2(t,s)$ in the form
\begin{equation} \label{sigma}
k_{2}(t,s)=n_{1}(t,s)+in_{2}(t,s),
\end{equation}
where
\[n_{1}(t,s):=\frac{1+st}{(s-t)(s^2+1)}\times\]
\[\times \frac{\sigma_{1}'(S)\Bigl(d_{2}(S,T)-\sigma_{2}'(S)\Bigr)
+\sigma_{2}'(S)\Bigl(\sigma_{1}'(S)-d_{1}(S,T)\Bigr)}{2\left(d(S,T)\right)^{2}}\,,\]
\[ n_{2}(t,s):=\frac{\sigma_{1}'(S)\Bigl(d_{2}(S,T)-\sigma_{2}'(S)\Bigr)
+\sigma_{2}'(S)\Bigl(\sigma_{1}'(S)-d_{1}(S,T)\Bigr)}{2\left(d(S,T)\right)^{2}(s^2+1)}\,.\]

Now, estimates of the type (\ref{ner_K1'})  for $I[g,n_{1}](t)$,
$I[g,n_{2}](t)$ are established similarly to analogous estimates
for $I[g,m_1](t)$, $I[g,m_2](t)$, respectively. The lemma is
proved.

Consider the notations
\[k_{1}(\infty,s):=-\frac{s}{s^2+1}\,\frac{\sigma'(S)-d(S,1)}{d(S,1)}+i\,\frac{\sigma'(S)}{(s^2+1)\,d(S,1)}=:\]
\[=:m_{1}(\infty,s)+im_{2}(\infty,s),\]
\[k_{2}(\infty,s):=-\frac{s}{2(s^2+1)
\left(d(S,1)\right)^{2}}\times\]
\[\times \Bigl(\sigma'(S)\bigl(d_{2}(S,1)-\sigma_{2}'(S)\bigr) +\sigma_{2}'(S)
\bigl(\sigma'(S)-d(S,1)\bigr) \Bigr)+\]
\[+i\,\frac{1}{2(s^2+1)\,\left(d(S,1)\right)^{2}} \Bigl(\sigma'(S)\, d_{2}(S,1)
-\sigma_{2}'(S)\,d(S,1)\Bigr)=:\]
\[=:n_{1}(\infty,s)+in_{2}(\infty,s).\]

 \begin{lem} \label{dini_K_infty}
Let $g\in C(\,\overline{\mathbb{R}}\,)$ and the conformal mapping
$\sigma(T)$ have the nonvanishing continuous contour derivative
$\sigma'(T)$ on the circle $\Gamma$, and its modulus of continuity
satisfy the condition \eqref{Dini_circle}. Then for
$0<\varepsilon<1/4$ and $t\in\mathbb{R}$ such that
$|\,t\,|>1/\varepsilon$ the following estimates are true:
 \[\left|\,\int\limits_{-\infty}^{\infty}g(s)\,k_{j}(t,s)\,ds-\int\limits_{-\infty}^{\infty}g(s)\,k_{j}
 (\infty,s)\,ds\,\right| \le\]
\begin{equation} \label{ner_K-m'_infty}
 \le c\,\|g\|_{C(\overline{\mathbb{R}}\,)}
 \,\varepsilon\int\limits_{0}^{2}\frac{\omega_{\Gamma}(\sigma', \eta)}{\eta(\eta+\varepsilon)}\,
 d\eta\,,\qquad j=1,2\,,
      \end{equation}
where the constant\, $c$ does not depend on $t$ and $\varepsilon$.
\end{lem}
{\bf Proof.} Consider the case  $t>1/\varepsilon$ (the case
$t<-1/\varepsilon$ is considered by analogy).

In order to prove the estimate  (\ref{ner_K-m'_infty}) for $j=1$,
we shall use the expression \eqref{K1-int} of the integral
$I[g,k_1](t)$ and obtain estimates of the type
\eqref{ner_K-m'_infty} for each of the integrals  $I[g,m_1](t)$,
$I[g,m_2](t)$.

For $I[g,m_1](t)$ we have
\[
 \left|\,\int\limits_{-\infty}^{\infty}g(s)\,m_{1}(t,s)\,ds-\int\limits_{-\infty}^{\infty}g(s)\,m_{1}(\infty,s)\,ds\,\right| \le \]
\[\le
\left(\int\limits_{-\infty}^{-t/2}+\int\limits_{3t/2}^{\infty}\right)
\left|g(s)\right|\Bigl(\left|m_{1}(t,s)\right|+\left|m_{1}(\infty,s)\right|\Bigr)
\, ds+ \]
\[+
\int\limits_{t/2}^{3t/2}
\left|g(s)\right|\Bigl(\left|m_{1}(t,s)\right|+\left|m_{1}(\infty,s)\right|\Bigr)
\, ds+
\]
\[+\int\limits_{-t/2}^{t/2}
\left|g(s)\right|\bigl|m_{1}(t,s)-m_{1}(\infty,s)\bigr|\,
ds=:I_1+I_2+I_3.\]

The integrals $I_1$ and $I_2$ are estimated with using
(\ref{sigma'_1}):
\[I_{1}\le c \,\|g\|_{C(\overline{\mathbb{R}}\,)}
 \, \left( \Biggl(\int\limits_{-\infty}^{-t/2}+\int\limits_{3t/2}^{\infty}\Biggr)\omega_{\Gamma}(\sigma',
|S-T|) \frac{t\,ds}{s^2+1}\right.+\]
\[+\left.\Biggl(\int\limits_{-\infty}^{-t/2}+\int\limits_{3t/2}^{\infty}\Biggr)\omega_{\Gamma}(\sigma',
|S-1|)\frac{ds}{\sqrt{s^2+1}} \right)\le \]
\[\le c
\,\|g\|_{C(\overline{\mathbb{R}}\,)} \,\left(
\omega_{\Gamma}(\sigma',6\varepsilon)
\Biggl(\int\limits_{-\infty}^{-t/2}+\int\limits_{3t/2}^{\infty}\Biggr)
\frac{t\,ds}{s^2}\right.+\]
\[+\left. \Biggl(\int\limits_{-\infty}^{-t/2}+\int\limits_{3t/2}^{\infty}\Biggr)\frac{\omega_{\Gamma}(\sigma',
|S-1|)}{|S-1|}\, \frac{|\,s|\,ds}{(s^2+1)^{3/2}}\right) \le\]
\[
\le c \,\|g\|_{C(\overline{\mathbb{R}}\,)}
 \,\left(
\omega_{\Gamma}(\sigma',6\varepsilon)+\int\limits_{0}^{4\varepsilon}\frac{\omega_{\Gamma}(\sigma',
\eta)}{\eta}\, d\eta\right) \le
c\,\|g\|_{C(\overline{\mathbb{R}}\,)}
 \,\varepsilon\int\limits_{0}^{2}\frac{\omega_{\Gamma}(\sigma', \eta)}{\eta(\eta+\varepsilon)}\,
 d\eta\,,\]
\[I_{2}\le c\,\|g\|_{C(\overline{\mathbb{R}}\,)}\, \int\limits_{t/2}^{3t/2}
\left( \frac{\omega_{\Gamma}(\sigma', |S-T|)}{|S-T|}
\frac{1+|s||\,t|}{\sqrt{s^2+1}\sqrt{t^2+1}}\right.+\]
\[+\biggl.\omega_{\Gamma}(\sigma', |S-1|)\,s
\biggr)\frac{ds}{s^2+1} \le\]
\[\le c\,\|g\|_{C(\overline{\mathbb{R}}\,)}\,\left( \int\limits_{0}^{2\varepsilon}
\frac{\omega_{\Gamma}(\sigma', \eta)}{\eta}\, d\eta+
 \omega_{\Gamma}(\sigma',4 \varepsilon)\right) \le c\,\|g\|_{C(\overline{\mathbb{R}}\,)}
 \,\varepsilon\int\limits_{0}^{2}\frac{\omega_{\Gamma}(\sigma', \eta)}{\eta(\eta+\varepsilon)}\,
 d\eta\,. \]
Here and below in the proof, by\, $c$\, we denote constants whose
values are independent of  $t$ and $\varepsilon$, but, generally
speaking, may be different even within a single chain of
inequalities.

For estimating the integral $I_3$ we use the equality
\[m_{1}(t,s)-m_{1}(\infty,s)=\left(\frac{1}{s-t}-\frac{s}{s^2+1}\right)\frac{\sigma'(S)-d(S,T)}{d(S,T)}
+\]
\[+\frac{s}{s^2+1}\,\frac{\sigma'(S)-d(S,1)}{d(S,1)}=
\]
\[
=\frac{s}{s^2+1}\,\sigma'(S)\,\frac{d(S,T)-d(S,1)}{d(S,1)\,d(S,T)}+
\frac{1}{s-t}\frac{\sigma'(S)-d(S,T)}{d(S,T)}\,,\]
 and the inequalities  \eqref{prir_d_1_k}, (\ref{sigma'_1}) and
properties of a modulus of continuity  (cf., e.g.,
\cite[p.~176]{Dzjadyk}). Thus, we obtain
\[I_{3} \le c\,\|g\|_{C(\overline{\mathbb{R}}\,)}\,
\left( |T-1|
\int\limits_{-t/2}^{t/2}\frac{\omega_{\Gamma}(\sigma',
|S-1|)}{|S-1|}\, \frac{|s|}{s^2+1}\, ds\right.+\]
\[+\left. \int \limits_{-t/2}^{t/2} \frac{\omega_{\Gamma}(\sigma',
|S-T|)}{|S-T|}\, \frac{ds}{\sqrt{s^2+1}\sqrt{t^2+1}} \right)\le \]
\[\le c\,\|g\|_{C(\overline{\mathbb{R}}\,)}\, |T-1| \int\limits_{-t/2}^{t/2}\frac{\omega_{\Gamma}(\sigma',
|S-1|)}{|S-1|^2}\, \frac{|\,s|\,ds}{(s^2+1)^{3/2}}  \le\]
\[\le c\,\|g\|_{C(\overline{\mathbb{R}}\,)}
 \,\varepsilon\int\limits_{\varepsilon}^{2}\frac{\omega_{\Gamma}(\sigma', \eta)}{\eta^2}\,
 d\eta \le
 c\,\|g\|_{C(\overline{\mathbb{R}}\,)}
 \,\varepsilon\int\limits_{0}^{2}\frac{\omega_{\Gamma}(\sigma', \eta)}{\eta(\eta+\varepsilon)}\,
 d\eta\,.\]

An estimate of the type \eqref{ner_K-m'_infty} for $I[g,m_1](t)$
follows from the obtained inequalities. An estimate of the type
\eqref{ner_K-m'_infty} for $I[g,m_2](t)$ is similarly established.
Thus, the inequality (\ref{ner_K-m'_infty}) is proved for $j=1$.

To prove the estimate (\ref{ner_K-m'_infty}) for $j=2$, we use the
representation \eqref{sigma} of the function $k_2(t,s)$ and obtain
estimates of the type \eqref{ner_K-m'_infty} for each of the
integrals $I[g,n_1](t)$, $I[g,n_2](t)$ similarly to the estimation
of $I[g,m_1](t)$. The lemma is proved.

The next statement follows obviously from Lemmas  \ref{dini_K1}
and \ref{dini_K_infty}.

\begin{thm} \label{copm_op-ra}
Let the conformal mapping $\sigma(T)$ have the nonvanishing
continuous contour derivative $\sigma'(T)$ on the circle $\Gamma$,
and its modulus of continuity satisfy the condition
\eqref{Dini_circle}. Let the function $k(t,s)$ be defined by the
relation \eqref{k}, where the functions
 $k_1(t,s)$, $k_2(t,s)$ are
 defined by the equalities
\eqref{k1}, \eqref{k2}, respectively. Then the operator
\[\EuScript{I}[g]:=\int\limits_{-\infty}^{\infty}g(s)\,k(t,s)\,ds\]
is compact in the space  $C(\,\overline{\mathbb{R}}\,)$.
\end{thm}
\section{Equivalence conditions of the  (1-3)-problem to a
system of Fredholm integral equations} \label{equiv-13}
 From the above, one can see that finding a solution of the (1-3)-problem in the form
\eqref{int_t_Cauchy-1-3} with functions  $\varphi_1 :
\partial D_{\zeta}\rightarrow \mathbb{R}$, $\varphi_3 : \partial
D_{\zeta}\rightarrow \mathbb{R}$ satisfying conditions of the type
\eqref{usl_Dini} is reduced to solving the system of integral
equations \eqref{syst-int-ur}. Under conditions of Theorem
\ref{copm_op-ra}, integral operators, which are generated by the
left parts of system \eqref{syst-int-ur}, are compact in the space
$C(\,\overline{\mathbb{R}}\,)$, i.e. the system
\eqref{syst-int-ur} is a system of Fredholm integral equations.

For any function   $g\in C(\,\overline{\mathbb{R}}\,)$ we use the
local centered (with respect to the infinitely remote point)
modulus of continuity
\[\omega_{\mathbb{R}, \infty}(g,\varepsilon):=\sup_{t\in
\mathbb{R}\, :\, |t|\ge 1/\varepsilon} |g(t)-g(\infty)|\,.\]

Let ${\cal D}(\overline{\mathbb{R}})$ denote the class of
functions $g\in C(\overline{\mathbb{R}}\,)$ whose moduli of
continuity satisfy the Dini conditions
 \begin{equation} \label{Dini_both}
\int\limits_{0}^{1}\frac{\omega_{\mathbb{R}}(g,\eta)}{\eta}\,d\eta<\infty,\qquad
\int\limits_{0}^{1}\frac{\omega_{\mathbb{R},\infty}(g,\eta)}{\eta}\,d\eta<\infty\,.
\end{equation}

Since the functions $\varphi_1 \colon \partial
D_{\zeta}\longrightarrow \mathbb{R}$ and $\varphi_3 \colon
\partial D_{\zeta}\longrightarrow \mathbb{R}$ in the expression \eqref{int_t_Cauchy-1-3} of a solution of
the (1-3)-problem have to satisfy conditions of the type
\eqref{usl_Dini}, it is necessary to require that corresponding
functions  $g_{1}$, $g_3$ satisfying the system
\eqref{syst-int-ur} should belong to the class ${\cal
D}(\overline{\mathbb{R}})$. In the next theorem we state a
condition on the conformal mapping $\sigma(T)$, under which all
solutions of the system \eqref{syst-int-ur} satisfy the mentioned
requirement.

 \begin{thm} \label{syst-CtoD}  Let the functions  $u_1 \colon
\partial D_{\zeta}\longrightarrow \mathbb{R}$, $u_3 \colon \partial
D_{\zeta}\longrightarrow \mathbb{R}$ satisfy conditions of the
type \eqref{usl_Dini}. Let the conformal mapping $\sigma(T)$ have
the nonvanishing continuous contour derivative $\sigma'(T)$ on the
circle $\Gamma$, and its modulus of continuity satisfy the
condition
\begin{equation}\label{usl_ln}
 \int\limits_{0}^{2}\frac{\omega_{\Gamma}(\sigma', \eta) }{\eta}\,\ln\frac{3}{\eta} \, d\eta < \infty.
 \end{equation}
Then all continuous functions $g_1, g_3$ satisfying the system of
Fredholm integral equations \eqref{syst-int-ur} belong to the
class  ${\cal D}(\overline{\mathbb{R}})$, and the corresponding
functions $\varphi_1$, $\varphi_3$ in \eqref{int_t_Cauchy-1-3}
satisfy conditions of the type \eqref{usl_Dini}.
\end{thm}
{\bf Proof.} Let us rewrite the system  \eqref{syst-int-ur} in the
matrix form:
\begin{equation}\label{vekt-syst}
\left(\begin{array}{l} g_1(t) \\[2mm] g_3(t)  \end{array}\right)=
\left(\begin{array}{l} 2\widetilde{u}_{1}(t) \\[2mm]
2\widetilde{u}_{3}(t)  \end{array}\right)-\left(\begin{array}{l}
U_{1} \left[ \frac{1}{\pi i}\,I[g,k](t)\right] \\[2mm] U_{3} \left[
\frac{1}{\pi i}\,I[g,k](t)\right]
\end{array}\right) \qquad
\forall\,t\in\mathbb{R}\,,
 \end{equation}
where \,\, $I[g,k](t):=\int\limits_{-\infty}^{\infty}g(s)k(t,s)\,
ds$ \, and the function $g$ is defined by the equality
\eqref{fun_g}.

Inasmuch as the derivative  $\sigma'(T)$ is continuous on $\Gamma$
and the functions  $u_1$, $u_3$ satisfy conditions of the type
\eqref{usl_Dini}, 
the right-hand sides $\widetilde u_1$, $\widetilde u_3$ of the
equations \eqref{syst-int-ur} belong to the class  ${\cal
D}(\overline{\mathbb{R}})$. With using Lemmas~\ref{dini_K1} and
\ref{dini_K_infty}, it is easy to establish that moduli of
continuity of the function $I[g,k](t)$ satisfy conditions of the
type \eqref{Dini_both} due to the condition \eqref{usl_ln}.
Therefore, in view of (\ref{vekt-syst}) the functions $g_1$, $g_3$
belong to the class ${\cal D}(\overline{\mathbb{R}})$ also.
Finally, by Lemma  3.3  in \cite{Mel_Pl_mon}, we conclude that the
functions  $\varphi_1$, $\varphi_3$ satisfy conditions of the type
\eqref{usl_Dini}. The theorem is proved.

Let us make some remarks concerning the representation of a
solution of the (1-3)-problem by the formula
\eqref{int_t_Cauchy-1-3}.

Suppose that a solution $\Phi$ of the (1-3)-problem is
continuously extended to the boundary $\partial D_{\zeta}$. By the
(2-4)-problem {\em conjugated} with the (1-3)-problem or, briefly,
the {\em conjugated} (2-4)-problem  we shall call a problem on
finding a continuous function
 $\Phi_{\ast}\colon \mu \setminus
D_{\zeta}\longrightarrow\mathbb{B}$ which is monogenic in the
domain  $\mu \setminus\overline{D_{\zeta}}$ and vanishes at the
infinity with the following boundary conditions:
\[U_{2}[\Phi_{\ast}(\zeta)]=U_{2}[\Phi(\zeta)], \quad
U_{4}[\Phi_{\ast}(\zeta)]=U_{4}[\Phi(\zeta)] \qquad \forall \zeta
\in \partial D_{\zeta}\,.\]

Note that if a solution of the (1-3)-problem is expressed by the
formula \eqref{int_t_Cauchy-1-3}, where the functions $\varphi_1$,
$\varphi_3$ satisfy conditions of the type \eqref{usl_Dini}, then
by Theorem \ref{form_Sokh_D}, the integral \eqref{int_t_Cauchy}
with $\varphi(\tau)=\varphi_1(\tau)e_1+\varphi_3(\tau)e_2$, being
a solution of the (1-3)-problem, can be continuously extended to
the boundary $\partial D_{\zeta}$ from each of the domains
$D_{\zeta}$, $\mu\setminus\overline{D_{\zeta}}$ and is also a
solution of the conjugated (2-4)-problem due to the formulas
\eqref{int_t_Cauchy+}.

Furthermore, by virtue of assumptions that a solution  $\Phi$ of
the (1-3)-problem is continuously extended to the boundary
$\partial D_{\zeta}$ and the conjugated (2-4)-problem has a
solution $\Phi_{\ast}$, it follows that the function  $\Phi$ is
represented in the form \eqref{int_t_Cauchy-1-3}. Indeed, by the
Cauchy integral formula and the Cauchy theorem for monogenic
functions in the biharmonic plane (cf., e.g, Theorem~3.2 in
\cite{Conrem_13}), the following equalities are true:
\[ \Phi(\zeta)=\frac{1}{2 \pi i}\int\limits_{\partial D_{\zeta}}\Phi(\tau)(\tau-\zeta)^{-1}\, d\tau,
\quad 0 = \frac{1}{2 \pi i}\int\limits_{\partial
D_{\zeta}}\Phi_{\ast}(\tau)(\tau-\zeta)^{-1}\, d\tau   \quad
\forall \zeta \in D_{\zeta}, \]
 that implies the equality \eqref{int_t_Cauchy-1-3} with
\[\varphi_{j}(\tau)=U_{j}\left[\Phi(\tau)\right]-U_{j}\left[\Phi_{\ast}(\tau)\right]=u_j(\tau)-U_{j}\left[\Phi_{\ast}(\tau)\right],\qquad
j\in\{1,3\}.\]

The made remarks can be amplified by the following theorem.

 \begin{thm} \label{osn-teor-a}  Let the functions  $u_1 :
\partial D_{\zeta}\longrightarrow \mathbb{R}$, $u_3 : \partial
D_{\zeta}\longrightarrow \mathbb{R}$ satisfy conditions of the
type \eqref{usl_Dini}. Let the conformal mapping $\sigma(T)$ have
the nonvanishing continuous contour derivative $\sigma'(T)$ on the
circle $\Gamma$, and its modulus of continuity satisfy the
condition \eqref{usl_ln}. Then the following assertions are
equivalent\emph{:}
\begin{enumerate}
\item[\emph{(I)}] the system of Fredholm  integral equations  \eqref{syst-int-ur}
is solvable in the space  $C(\overline{\mathbb{R}})$;
\item[\emph{(II)}] there exists a solution of the \em (1-3)-\em problem of the form
\eqref{int_t_Cauchy-1-3}, where the functions  $\varphi_1$,
$\varphi_3$ satisfy conditions of the type
\eqref{usl_Dini}\emph{;}
\item[\emph{(III)}]
   a solution $\Phi$ of
   the \em (1-3)-\em problem is continuously extended to the boundary $\partial
   D_{\zeta}$. For this function $\Phi$,
   the conjugated \em (2-4)-\em problem is solvable and moduli of continuity of components $U_{1}\left[\Phi_{\ast}\right]$, $U_{3}\left[\Phi_{\ast}\right]$
   of its solution  $\Phi_{\ast}$
   satisfy conditions of the type \eqref{usl_Dini}.
\end{enumerate}
\end{thm}
{\bf Proof.} Continuing the reasonings adduced before
Theorem~\ref{osn-teor-a}, we conclude that in the case where the
functions $u_1$, $u_3$ satisfy conditions of the type
\eqref{usl_Dini}, the functions $\varphi_1$, $\varphi_3$ satisfy
the same conditions if and only if moduli of continuity of
components $U_{1}\left[\Phi_{\ast}\right]$,
$U_{3}\left[\Phi_{\ast}\right]$ of a solution $\Phi_{\ast}$ of the
conjugated  (2-4)-problem satisfy conditions of the type
 \eqref{usl_Dini}. Thus, the assertions  (II) and (III) are equivalent.
The equivalence of assertions (I) and (II) is a consequence of
Theorem~\ref{syst-CtoD}. The theorem is proved.

Rewrite the integral equations of the system  \eqref{syst-int-ur}
in expanded form:
\[\frac{1}{2}\,g_{1}(t)+\frac{1}{2 \pi}\,
\mathrm{Im}\int\limits_{-\infty}^{\infty}g_{1}(s)\frac{\tau'(s)}{\tau(s)-\tau(t)}\,ds + 
\frac{1}{2 \pi}\,
\mathrm{Re}\int\limits_{-\infty}^{\infty}g_{1}(s)\frac{\tau_{1}'(s)(\tau_{2}(s)-\tau_{2}(t))}
{(\tau(s)-\tau(t))^{2}}\,ds-
\]
\[
-\frac{1}{2\pi}\,\mathrm{Re}\int\limits_{-\infty}^{\infty}
g_{1}(s)
\frac{\tau_{2}'(s)(\tau_{1}(s)-\tau_{1}(t))}{(\tau(s)-\tau(t))^{2}}\,ds-\]
\[-
\frac{1}{2 \pi}\,
\mathrm{Im}\int\limits_{-\infty}^{\infty}g_{3}(s)
\frac{\tau_{1}'(s)(\tau_{2}(s)-\tau_{2}(t))}
{(\tau(s)-\tau(t))^{2}}\,ds+\]
\[+
\frac{1}{2\pi}\, \mathrm{Im}\int\limits_{-\infty}^{\infty}g_{3}(s)
\frac{\tau_{2}'(s)(\tau_{1}(s)-\tau_{1}(t))}{(\tau(s)-\tau(t))^{2}}\,ds=\widetilde{u}_{1}(t),\]
\[\frac{1}{2}\,g_{3}(t)+\frac{1}{2 \pi}\, \mathrm{Im}\int\limits_{-\infty}^{\infty}g_{3}(s)
\frac{\tau'(s)}{\tau(s)-\tau(t)}\,ds-
\frac{1}{2 \pi}\,
\mathrm{Re}\int\limits_{-\infty}^{\infty}g_{3}(s)\frac{\tau_{1}'(s)(\tau_{2}(s)-\tau_{2}(t))}
{(\tau(s)-\tau(t))^{2}}\,ds+
\]
\[+\frac{1}{2\pi}\,\mathrm{Re}\int\limits_{-\infty}^{\infty}g_{3}(s)
\frac{\tau_{2}'(s)(\tau_{1}(s)-\tau_{1}(t))}{(\tau(s)-\tau(t))^{2}}\,ds-\]
\[- \frac{1}{2 \pi}\,
\mathrm{Im}\int\limits_{-\infty}^{\infty}g_{1}(s)\frac{\tau_{1}'(s)(\tau_{2}(s)-\tau_{2}(t))}
{(\tau(s)-\tau(t))^{2}}\,ds+\] 
\[+\frac{1}{2\pi}\,
\mathrm{Im}\int\limits_{-\infty}^{\infty}g_{1}(s)
\frac{\tau_{2}'(s)(\tau_{1}(s)-\tau_{1}(t))}{(\tau(s)-\tau(t))^{2}}\,ds=\widetilde{u}_{3}(t)
\qquad \forall\,t\in\mathbb{R}.
\]

Then the homogeneous system transposed with the system
\eqref{syst-int-ur} is of the form:
\[h_{1}(t)-\frac{1}{\pi}\, \mathrm{Im}\left(\tau'(t)\int\limits_{-\infty}^{\infty}
\frac{h_{1}(s)\,ds }{\tau(s)-\tau(t)}\right)-
\frac{\tau_{1}'(t)}{\pi}\,
\mathrm{Re}\int\limits_{-\infty}^{\infty}h_{1}(s)\frac{(\tau_{2}(s)-\tau_{2}(t))}
{(\tau(s)-\tau(t))^{2}}\,ds+
\]
\[
+\frac{\tau_{2}'(t)}{\pi}\,\mathrm{Re}\int\limits_{-\infty}^{\infty}h_{1}(s)
\frac{(\tau_{1}(s)-\tau_{1}(t))}{(\tau(s)-\tau(t))^{2}}\,ds+
\frac{\tau_{1}'(t)}{\pi}\,
\mathrm{Im}\int\limits_{-\infty}^{\infty}h_{3}(s)\frac{(\tau_{2}(s)-\tau_{2}(t))}
{(\tau(s)-\tau(t))^{2}}\,ds-\]
\begin{equation} \label{h1}
-\frac{\tau_{2}'(t)}{\pi}\,
\mathrm{Im}\int\limits_{-\infty}^{\infty}h_{3}(s)
\frac{(\tau_{1}(s)-\tau_{1}(t))}{(\tau(s)-\tau(t))^{2}}\,ds=0,
\end{equation}
\[h_{3}(t)-\frac{1}{\pi}\, \mathrm{Im}\left(\tau'(t)\int\limits_{-\infty}^{\infty}
\frac{h_{3}(s)\,ds}{\tau(s)-\tau(t)}\right)
+\frac{\tau_{1}'(t)}{\pi}\,
\mathrm{Re}\int\limits_{-\infty}^{\infty}h_{3}(s)\frac{(\tau_{2}(s)-\tau_{2}(t))}
{(\tau(s)-\tau(t))^{2}}\,ds-
\]
\[-\frac{\tau_{2}'(t)}{\pi}\,\mathrm{Re}\int\limits_{-\infty}^{\infty}h_{3}(s)
\frac{(\tau_{1}(s)-\tau_{1}(t))}{(\tau(s)-\tau(t))^{2}}\,ds +
\frac{\tau_{1}'(t)}{\pi}\,
\mathrm{Im}\int\limits_{-\infty}^{\infty}h_{1}(s)\frac{(\tau_{2}(s)-\tau_{2}(t))}
{(\tau(s)-\tau(t))^{2}}\,ds-\]
\begin{equation} \label{h3}
-\frac{\tau_{2}'(t)}{\pi}\,
\mathrm{Im}\int\limits_{-\infty}^{\infty}h_{1}(s)
\frac{(\tau_{1}(s)-\tau_{1}(t))}{(\tau(s)-\tau(t))^{2}}\,ds=0\qquad
\forall\,t\in\mathbb{R}\,.
\end{equation}

\begin{lem} \label{sol-souzn_syst} The pair of functions  $(h_1, h_2):=(\tau_{1}',\tau_{2}')$
satisfies the system of integral equations \eqref{h1}, \eqref{h3}.
\end{lem}
{\bf Proof.}  Substituting $h_{1}=\tau_{1}'$, $h_{2}=\tau_{2}'$
into the equation (\ref{h1}) and applying the conformal mapping
$z=\tau(t)$ of the upper half-plane  $\{t\in \mathbb{C}:
\mathrm{Im}\,t>0\}$ onto the domain $D_{z}$ of complex plane, we
pass in  (\ref{h1}) to integrating along the boundary $\partial
D_{z}$ of domain $D_{z}$. In such a way, using the denotations
$v_1:=\tau_{1}(s)$, $v_2:=\tau_{2}(s)$, $v:=\tau(s)\equiv v_1+i
v_2$, $x:=\tau_{1}(t)$, $y:=\tau_{2}(t)$, $z \equiv x+i y $, we
obtain
\[\tau_{1}'(t)-\frac{1}{\pi}\,\mathrm{Im}\left(\tau'(t)\int\limits_{\partial D_z}\frac{dv_1}{v-z}
\right)-\]
\[- \frac{\tau_{1}'(t)}{\pi}\,\mathrm{Re}\int\limits_{\partial
D_z}\frac{v_2-y}{(v-z)^2}\,dv_1+
\frac{\tau_{2}'(t)}{\pi}\,\mathrm{Re}\int\limits_{\partial
D_z}\frac{v_1-x}{(v-z)^2}\,dv_1+\]
\begin{equation}\label{souzn-1}
+\frac{\tau_{1}'(t)}{\pi}\,\mathrm{Im}\int\limits_{\partial
D_z}\frac{v_2-y}{(v-z)^2}\,dv_2-
\frac{\tau_{2}'(t)}{\pi}\,\mathrm{Im}\int\limits_{\partial
D_z}\frac{v_1-x}{(v-z)^2}\,dv_2=0.
\end{equation}
Using polar coordinates $(r, \theta)$ with the relation
$v-z=:r\exp\{i\theta\}$, it is easy to show that the equality
\eqref{souzn-1} becomes identical.

In a similar way, we conclude that the functions
$h_{1}=\tau_{1}'$, $h_{2}=\tau_{2}'$ satisfy the equation
(\ref{h3}). The lemma is proved.

It follows from Lemma~\ref{sol-souzn_syst} that the condition 
\[\int\limits_{-\infty}^{\infty} \Bigl(\widetilde{u}_1(s)\tau_1'(s)+\widetilde{u}_3(s)\tau_2'(s)\Bigr)\,ds=0 \]
is necessary for the solvability of the system of integral
equations \eqref{syst-int-ur}. Passing in this equality to the
integration along the boundary $\partial D_{\zeta}$, we get the
equivalent condition
\begin{equation} \label{Razhr}
\int\limits_{\partial D_{\zeta}} u_1(xe_1+ye_2) \, d x+
u_3(xe_1+ye_2) \, dy =0
\end{equation}
for given functions $u_1$, $u_3$.

It is proved in \cite{IJPAM_13}  that the condition  \eqref{Razhr}
is also sufficient for the solvability of the (1-3)-problem for a
disk. The next theorem contains assumptions, under which the
condition  \eqref{Razhr} is necessary and sufficient for the
solvability of the system of integral equations
\eqref{syst-int-ur} and, therefore, for the existence of a
solution of the (1-3)-problem in the form
\eqref{int_t_Cauchy-1-3}.

\begin{thm} \label{edinstv-suF} Assume that the conformal mapping $\sigma(T)$
has the nonvanishing continuous contour derivative $\sigma'(T)$ on
the circle $\Gamma$, and its modulus of continuity satisfies the
condition \eqref{usl_ln}. Also, assume that all solutions $g_1,
g_3$ of the homogeneous system of equations \eqref{syst-int-ur}
\em (\em with $\widetilde{u}_{j} \equiv 0$ for $j\in\{1,3\}$\em )
\em are differentiable on $\mathbb{R}$, and the integral
\[
\frac{1}{2\pi i}\int\limits_{\partial
D_\zeta}\varphi'(\tau)(\tau-\zeta)^{-1}\,d\tau\]
 is bounded in
both domains  $D_{\zeta}$ and $\mu\setminus\overline{D_{\zeta}}$
\em ; \em here $\varphi'$ is the contour derivative of the
function
 $\varphi(\tau):=g_{1}(s)e_1+g_{3}(s)e_2$, where  $\tau=\widetilde{\tau}(s)$ for all
$s\in\mathbb{R}$.
 Then the following assertions are true:
\begin{enumerate}
\item[\emph{1)}] the number of linearly independent solutions of the homogeneous system of equations  \eqref{syst-int-ur} is
equal to $1$;
\item[\emph{2)}]
the non-homogeneous system of equations \eqref{syst-int-ur} is
solvable if and only if the condition \eqref{Razhr} is satisfied.
\end{enumerate}
\end{thm}
{\bf Proof.} It follows from Lemma~\ref{sol-souzn_syst} that the
transposed system of equations  (\ref{h1}), (\ref{h3}) has at
least one nontrivial solution. Therefore, by the Fredholm theory,
the system of homogeneous equations  \eqref{syst-int-ur} has at
least one nontrivial solution $g_1=g_1^{0}$, $g_3=g_{3}^{0}$.
Consider the function
\begin{equation}\label{def-phi0}
\varphi_0(\tau):=g_{1}^{0}(s)e_1+g_{3}^{0}(s)e_2\,, \quad
\tau=\widetilde{\tau}(s)\,, \quad \forall\,s\in\mathbb{R}\,
\end{equation}
corresponding to this solution.

Then the function, which is defined for all $\zeta\in D_{\zeta}$
by the formula
\begin{equation}\label{phi0}
\Phi_{0}(\zeta):=\frac{1}{2\pi i}\int\limits_{\partial
D_\zeta}\varphi_0(\tau)(\tau-\zeta)^{-1}\,d\tau\,,
\end{equation}
is a solution of the homogeneous  (1-3)-problem (with
$u_1=u_3\equiv 0$). By virtue of differentiability of the
functions $g_1^{0}$, $g_{3}^{0}$ on $\mathbb{R}$, the function
$\varphi_0$ has the contour derivative $\varphi_0'(\tau)$ for all
$\tau\in\partial D_\zeta$ and the following equality holds:
\begin{equation}\label{phi0'}
\Phi_{0}'(\zeta):=\frac{1}{2\pi i}\int\limits_{\partial
D_\zeta}\varphi_0'(\tau)(\tau-\zeta)^{-1}\,d\tau \qquad \forall\,
\zeta\in D_{\zeta}\,.
\end{equation}

Let $\Phi_{1}$ be a monogenic in 
$D_{\zeta}$ function such that $\Phi_{1}'(\zeta)=\Phi_{0}(\zeta)$
for all $\zeta\in D_{\zeta}$. Then the function
$V(x,y):=U_1[\Phi_{1}(\zeta)]$ is a solution of the homogeneous
biharmonic problem \eqref{osn_big_pr} (with
$\omega_3=\omega_4\equiv 0$).

Taking into account the equalities
\[\Delta V(x,y):= \frac{\partial^{2} V(x,y)}{\partial
x^2}+\frac{\partial^{2} V(x,y)}{\partial y^2}=\]
\[=U_1\left[\Phi_{1}''(\zeta)(e_1^2+e_2^2)\right]=
U_1\left[\Phi_{0}'(\zeta)(e_1^2+e_2^2)\right]\] and boundedness of
the function \eqref{phi0'} in the domain  $D_\zeta$, we conclude
that the Laplacian
 $\Delta V(x,y)$ is bounded in the domain $D$.

By the uniqueness theorem (see \S~4 in \cite{Mikhlin_t_ed}),
solutions of the homogeneous biharmonic problem \eqref{osn_big_pr}
in the class of functions $V$, for which $\Delta V(x,y)$ is
bounded in $D$, are only constants: $V={\rm const}$. Therefore,
taking into account the fact that the functions \eqref{phi0} and
\eqref{osn-(1-3)} are equal, we obtain the equalities
$U_{1}[\Phi_{0}(\zeta)]=U_{3}[\Phi_{0}(\zeta)]\equiv 0$ in the
domain $D_\zeta$.

All monogenic functions $\Phi_{0}$ satisfying the condition
$U_{1}[\Phi_{0}(\zeta)]\equiv 0$ are described in
\cite[Lemma~3]{GrPl_umz-09}. Taking into account the identity
$U_{3}[\Phi_{0}(\zeta)]\equiv 0$ also, we obtain the equality
\begin{equation}\label{Phi_0}
\Phi_{0}(\zeta) = ik\zeta+i\left(n_1e_1+n_2e_2\right) \qquad
\forall\,\zeta \in \overline{D_{\zeta}}\,,
\end{equation}
where $k$, $n_1$ and $n_2$ are real numbers.

Let us show that the constant $k$ does not equal to zero in the
equality  (\ref{Phi_0}). Assume the contrary: $k=0$. Then it
follows from the Sokhotski-Plemelj formula \eqref{int_t_Cauchy+}
for the function \eqref{phi0} and the identity
$\Phi_{0}^+(\tau)\equiv i\left(n_1e_1+n_2e_2\right)$ that the
following equality holds:
\begin{equation}\label{Phi_0^-}
\Phi_{0}^-(\tau)=i\left(n_1e_1+n_2e_2\right)-\varphi_0(\tau)
\qquad \forall\,\tau\in \partial D_{\zeta}\,.
\end{equation}

Using the integral Cauchy formula for monogenic functions in the
biharmonic plane (cf., e.g., Theorem~3.2 in \cite{Conrem_13}) and
taking into account the equality  \eqref{Phi_0^-}, we obtain the
equalities
\[ \Phi_{0}(\zeta)=-\frac{1}{2 \pi i}\int\limits_{\partial D_{\zeta}}\Phi_{0}^-(\tau)(\tau-\zeta)^{-1}\,
d\tau=\]
\[=\frac{1}{2 \pi i}\int\limits_{\partial
D_{\zeta}}\varphi_0(\tau)(\tau-\zeta)^{-1}\,
d\tau-\frac{i\left(n_1e_1+n_2e_2\right)}{2 \pi
i}\int\limits_{\partial D_{\zeta}}(\tau-\zeta)^{-1}\, d\tau=\]
\[=\Phi_{0}(\zeta)-i\left(n_1e_1+n_2e_2\right) \qquad
\forall\,\zeta \in \mu\setminus\overline{D_{\zeta}}\,, \]
 a corollary of which is the equality\,\,  $n_1e_1+n_2e_2=0$. In
virtue of the uniqueness of decomposition of any element in the
algebra $\mathbb{B}$ with respect to the basis $\{e_1,e_2\}$, we
get $n_1=n_2=0$\,.

Thus, the equality $\Phi_{0}(\zeta)=0$ for all $\zeta \in
\overline{D_{\zeta}}$ is a corollary of the assumption $k=0$\,.
Furthermore, the function  \eqref{phi0}, which is considered for
all $\zeta\in\mu\setminus\overline{D_{\zeta}}$, is a solution of
the conjugated (2-4)-problem with boundary data:
\[U_{2}[\Phi_{0}^-(\zeta)]=U_{2}[\Phi_0^+(\zeta)]\equiv 0, \quad
U_{4}[\Phi_{0}^-(\zeta)]=U_{4}[\Phi_0^+(\zeta)]\equiv 0 \qquad
\forall \zeta \in \partial D_{\zeta}\,.\]

Using arguments analogous to those used in the proof of the
equality  \eqref{Phi_0}, we obtain the equality
\[\Phi_{0}(\zeta) = k_1\zeta+m_1e_1+m_2e_2 \qquad \forall\,\zeta \in
\mu\setminus\overline{D_{\zeta}},\]
 where $k_1$, $m_1$ and $m_2$ are real numbers. Therefore, as a
corollary of \eqref{Phi_0^-}, we get the equality
\begin{equation}\label{Phi_0-plot}
\varphi_0(\tau)=-\Phi_{0}^-(\tau)=-k_1\tau-m_1e_1-m_2e_2 \qquad
\forall\,\tau\in \partial D_{\zeta}\,.
\end{equation}

Substituting the expression  \eqref{Phi_0-plot} of the function
$\varphi_0$ in \eqref{phi0} and taking into account the equality
$\Phi_{0}(\zeta)=0$ for all $\zeta \in D_{\zeta}$, we obtain the
equality
\[0=-k_1\zeta-m_1e_1-m_2e_2 \qquad \forall\,\zeta \in
D_{\zeta},\]
 whence we get\,\, $k_1=m_1=m_2=0$. As a result, the
equality  \eqref{Phi_0-plot} is reduced to the identity
$\varphi_0(\tau)\equiv 0$ which contradicts to the non-triviality
of at least one of the functions $g_{1}^{0}$, $g_{3}^{0}$ in the
definition \eqref{def-phi0} of the function $\varphi_0$.

Hence, the assumption $k=0$ is not true. Therefore, the equality
(\ref{Phi_0}) holds with $k \ne 0$.

Let $g_1=g_1^{1}$, $g_3=g_{3}^{1}$ be a nontrivial solution of the
homogeneous system  \eqref{syst-int-ur}, and $g_1^{1}, g_{3}^{1}$
be different from  $g_1^{0}, g_{3}^{0}$. Then for the function
$\varphi^1(\tau):=g_{1}^1(s)e_1+g_{3}^1(s)e_2$, where
$\tau=\widetilde{\tau}(s)$ for all $s\in\mathbb{R}$, the following
equality of type (\ref{Phi_0}) holds:
\begin{equation}\label{Phi_0-new}
\frac{1}{2\pi i}\int\limits_{\partial
D_\zeta}\varphi^1(\tau)(\tau-\zeta)^{-1}\,d\tau=i\widetilde{k}\zeta+i\left(\widetilde{n}_1e_1+\widetilde{n}_2e_2\right)
\qquad \forall\,\zeta \in \overline{D_{\zeta}}\,,
\end{equation}
where $\widetilde{k}$, $\widetilde{n}_1$ and $\widetilde{n}_2$ are
real numbers.

Define a function  $\widetilde{\varphi}$ by the equality
\[\widetilde{\varphi}(\tau) =
\varphi^1(\tau)-\frac{\widetilde{k}}{k}\,\varphi_0(\tau) \qquad
\forall\,\tau\in \partial D_{\zeta}\,.\]

Taking into account the equalities (\ref{Phi_0}) and
(\ref{Phi_0-new}), we obtain
\[\frac{1}{2\pi i}\int\limits_{\partial D_\zeta}
\widetilde{\varphi}(\tau)(\tau-\zeta)^{-1}\,d\tau
=ie_1\left(\widetilde{n}_1-n_1\,\frac{\widetilde{k}}{k} \right)+
ie_2 \left(\widetilde{n}_2-n_2\, \frac{\widetilde{k}}{k}\right) \,
\forall \zeta \in \overline{D_{\zeta}}.\]

Further, in such a way as for the assumption $k=0$ in the equality
(\ref{Phi_0}), we obtain the identity
$\widetilde{\varphi}(\tau)\equiv 0$ which implies the equalities
$g_{m}^{1}(\tau)=\frac{\widetilde{k}}{k}\,g_{m}^{0}(\tau)$ for all
$\tau\in \partial D_{\zeta}$ and $m\in \{1,3\}$. The assertion  1)
is proved.

Due to the assertion  1), by the Fredholm theory, the number of
linearly independent solutions of the transposed system of
equations \eqref{h1}, \eqref{h3} is also equal to $1$. Such a
nontrivial solution is described in Lemma~\ref{sol-souzn_syst}.
Therefore, the condition \eqref{Razhr} is not only necessary but
also sufficient for the solvability of the system of integral
equations \eqref{syst-int-ur}. The theorem is proved.

\end{document}